# Almost Sure Convergence of Extreme Order Statistics

[a]Peng Zuoxiang    [a]Li Jiaona    [b]Saralees Nadarajah

[a]School of Mathematics and Statistics, Southwest University, Chongqing, 400715, China

[b]School of Mathematics, University of Manchester, Manchester, United Kingdom

**Abstract.** Let $M_n^{(k)}$ denote the $k$th largest maximum of a sample $(X_1, X_2, \cdots, X_n)$ from parent $X$ with continuous distribution. Assume there exist normalizing constants $a_n > 0$, $b_n \in \mathbb{R}$ and a nondegenerate distribution $G$ such that $a_n^{-1}(M_n^{(1)} - b_n) \xrightarrow{w} G$. Then for fixed $k \in \mathbb{N}$, the almost sure convergence of

$$\frac{1}{D_N} \sum_{n=k}^{N} d_n \mathbb{I}\{M_n^{(1)} \leq a_n x_1 + b_n, M_n^{(2)} \leq a_n x_2 + b_n, \cdots, M_n^{(k)} \leq a_n x_k + b_n\}$$

is derived if the positive weight sequence $(d_n)$ with $D_N = \sum_{n=1}^{N} d_n$ satisfies conditions provided by Hörmann.

**AMS Classification.** Primary 62F15; Secondary 60G70, 60F15.

**Keywords**: Almost sure convergence; Order statistics.

## 1 Introduction

Partial sum type almost sure (a.s.) central limit theorems were initiated by the papers of Brosamler (1988), Schatte (1988) and Lacey and Philipp (1990). For an i.i.d. sequence $\{X_n, n \geq 1\}$ with mean 0, variance 1 and partial sum $S_n = \sum_{k=1}^{n} X_k$, the simplest version of the a.s. central limit theorem says

$$\frac{1}{\log N} \sum_{n=1}^{N} \frac{1}{n} \mathbb{I}\{S_n \leq \sqrt{n} x\} \to \Phi(x), \quad a.s. \forall x \in \mathbb{R},$$

where $\mathbb{I}_A$ denotes the indicator function and $\Phi(x)$ is the standard normal distribution function. For unbounded functional a.s. central limit theorems, Ibragimov and Lifshits (1998) and Berkes et al. (1998) obtained the following a.s. central limit theorem under different restrictions on the continuous function $f$:

$$\frac{1}{\log N} \sum_{n=1}^{N} \frac{1}{n} f(S_n/\sqrt{n}) \to \int_{-\infty}^{\infty} f(x) d\Phi(x), \quad a.s.$$

The universal version of the a.s. central limit theorem discussed by Berkes and Csáki (2001) includes the case of the a.s. central limit theorem of extremes of i.i.d random sequences which was first studied by Fahrner and Stadtmüller (1998) and Cheng et al. (1998). Let $\{X_n, n \geq 1\}$ be an i.i.d. sequence, and let $M_n = \max_{1 \leq k \leq n} X_k$ denote the partial maximum. If there exist normalizing constants $a_n > 0$, $b_n \in \mathbb{R}$ and a nondegenerate distribution $G(x)$ such that

$$P(M_n \leq a_n x + b_n) \to G(x) \tag{1.1}$$

then

$$\frac{1}{\log N} \sum_{n=1}^{N} \frac{1}{n} \mathbb{I}\{M_n \leq a_n x + b_n\} \to G(x), \quad a.s. \forall x \in \mathbb{R}. \tag{1.2}$$

Fahrner (2000) extended (1.2) to unbounded continuous functions. The general result is

$$\frac{1}{\log N} \sum_{n=1}^{N} \frac{1}{n} f(a_n^{-1}(M_n - b_n)) \to \int_{-\infty}^{\infty} f(x) dG(x), \quad a.s.$$

See Theorem 1 of Fahrner (2000). Stadtmüller (2002) considered the a.s. central limit theorem of the $k$th maximum as $k = k_n$ satisfies $\log k_n = O((\log n)^{1-\varepsilon})$ for some $\varepsilon > 0$ or $(n - k_n)/n = p + O(1/\sqrt{n \log^\varepsilon n})$ for some $0 < p < 1$. Especially for fixed $k$ he showed

$$\frac{1}{\log N} \sum_{n=k}^{N} \frac{1}{n} \mathbb{I}\{M_n^{(k)} \leq a_n x + b_n\} \to G(x) \sum_{j=0}^{k-1} \frac{(-\log G(x))^j}{j!}, \quad a.s. \forall x \in \mathbb{R},$$

where $M_n^{(k)}$ denotes the $k$th maximum of $X_1, X_2, \cdots, X_n$ and $M_n := M_n^{(1)}$. Peng and Qi (2003) proved the a.s. central limit theorem of central order statistics, see also Hörmann (2005).

In this paper, we consider the following averages:

$$\frac{1}{D_N} \sum_{n=k}^{N} d_n \mathbb{I}\{M_n^{(1)} \leq a_n x_1 + b_n, M_n^{(2)} \leq a_n x_2 + b_n, \cdots, M_n^{(k)} \leq a_n x_k + b_n\}$$

provided the positive weights $d_n$, $n \geq 1$ satisfy the following conditions:

$$\liminf_{n \to \infty} n d_n > 0, \tag{1.3}$$

$$n^\alpha d_n \text{ is eventually nonincreasing for some } 0 < \alpha < 1, \tag{1.4}$$

and

$$\limsup_{n \to \infty} n d_n (\log D_n)^\rho / D_n < \infty \tag{1.5}$$

for some $\rho > 0$, where $D_n = \sum_{k=1}^{n} d_k$. Under conditions (1.3)–(1.5) it follows from the results in Hörmann (2006, 2007a, 2007b) that

$$\frac{1}{D_N} \sum_{n=1}^{N} d_n \mathbb{I}\{S_n \leq \sqrt{n}x\} \to \Phi(x), \quad a.s.,$$

and

$$\frac{1}{D_N} \sum_{n=1}^{N} d_n \mathbb{I}\{M_n \leq a_n x + b_n\} \to G(x), \quad a.s.$$

As discussed by Berkes and Csáki (2001) and Hömann (2005, 2006, 2007a, 2007b), the larger the $D_n$, the stronger the almost sure central limit theorem (ASCLT). If $d_n < 1/n$ such that $D_n \to \infty$, then the ASCLT holds. If $d_n = 1$, there is no ASCLT on the partial sum and partial maxima. Conditions (1.3), (1.4) and (1.5) tell us that there exists a large class of sequences $1/n < d_n < 1$ such that the ASCLT holds. For example, we may assume $D_n \to \infty$ with Karamata representation

$$D_n = \exp\left(\int_a^n \frac{\theta(u)}{u} du\right), \quad n > a,$$

where $\theta(x)$ is a slowly varying function such that

$$\liminf_{n\to\infty} \theta(n)D_n > 0, \quad \text{and} \quad \limsup_{n\to\infty} \theta(n)(\log D_n)^\rho < \infty$$

for some $\rho > 0$, which guarantees that (1.3), (1.4) and (1.5) hold. By the mean value theorem, we may choose $d_n \sim \theta(n)D_n/n$. This implies that $(d_n)$ is a regularly varying function with index $-1$. We mention the following examples:

(a) $D_n = (\log n)^\kappa$ with $d_n \sim \kappa(\log n)^{\kappa-1}/n$ for $\kappa > 1$, $\rho > 0$;

(b) $D_n = \exp((\log n)^\kappa)$ with $d_n \sim \kappa \exp((\log n)^\kappa)(\log n)^{\kappa-1}/n$ for $0 < \kappa < 1$, $0 < \rho < (1-\kappa)/\kappa$;

(c) $D_n = (\log n)^{1-\kappa} \exp((\log n)^\kappa)$ with $d_n \sim \kappa \exp((\log n)^\kappa)/n$ for $0 \le \kappa < 1/2$, $0 < \rho < 1/\kappa - 1$.

Throughout this paper we assume that $F(x)$, the univariate marginal distribution of $X_n$, $n \ge 1$ is continuous. This assumption implies that the order statistics are a.s. uniquely defined. Before providing the main results, recall the joint limit distribution of $\left(M_n^{(1)}, M_n^{(2)}, \cdots, M_n^{(k)}\right)$ for fixed $k$ if (1.1) holds. Define levels $u_n(x_j) = a_n x_j + b_n, j = 1, 2, \cdots, k, x_1 > x_2 > \cdots > x_k$ and define the point process $\chi_n$ of exceedances of levels $u_n(x_j), j = 1, 2, \cdots, k$ by i.i.d random variables $X_1, X_2, \cdots, X_n$. Then $\chi_n$ converges in distribution to a Poisson process on $(0, 1] \times \mathbb{R}$, for more details see Chapter 5 of Leadbetter et al. (1983), which states the joint limit distribution of $\left(M_n^{(1)}, M_n^{(2)}, \cdots, M_n^{(k)}\right)$ and

$$P\left(M_n^{(j)} \le a_n x + b_n\right) \to G(x) \sum_{i=0}^{j-1} \frac{(-\log G(x))^i}{i!} =: H_j(x), \quad j = 1, 2, \cdots, k \tag{1.6}$$

as $n \to \infty$. The joint limit distribution of $\left(M_n^{(1)}, M_n^{(2)}, \cdots, M_n^{(k)}\right)$ is so complicated that we express it as $H(x_1, x_2, \cdots, x_k)$ with the marginal distribution $H_j(x)$ defined in (1.6), $j = 1, 2, \cdots, k$, i.e.

$$\lim_{n\to\infty} P\left(M_n^{(1)} \le a_n x_1 + b_n, M_n^{(2)} \le a_n x_2 + b_n, \cdots, M_n^{(k)} \le a_n x_k + b_n\right)$$
$$= \begin{cases} H(x_1, x_2, \cdots, x_k), & x_1 > x_2 > \cdots > x_k; \\ 0, & \text{otherwise.} \end{cases}$$

## 2 Main Results

In this section, we provide the main results. The proofs are deferred to the next section.

**Theorem 1.** *Suppose (1.1) holds for an i.i.d. random sequence $(X_n, n \geq 1)$. Further assume (1.3)–(1.5) hold for positive weights $d_n, n \geq 1$. Then for fixed $k \in \mathbb{N}$ and real numbers $x_1 > x_2 > \cdots > x_k$, we have*

$$\frac{1}{D_N} \sum_{n=k}^{N} d_n \mathbb{I}\{M_n^{(1)} \leq u_n(x_1), M_n^{(2)} \leq u_n(x_2), \cdots, M_n^{(k)} \leq u_n(x_k)\}$$
$$\to H(x_1, x_2, \cdots, x_k), \quad a.s.,$$

*where $u_n(x_j), j = 1, 2, \cdots, k$ and $H(x_1, x_2, \cdots, x_k)$ are defined as before.*

**Corollary 1.** *Under the conditions of Theorem 1, for real numbers $x_{k_1} > x_{k_2} > \cdots > x_{k_l}$ with $1 \leq k_1 < k_2 < \cdots < k_l \leq k$, we have*

$$\frac{1}{D_N} \sum_{n=k_l}^{N} d_n \mathbb{I}\{M_n^{(k_1)} \leq u_n(x_{k_1}), M_n^{(k_2)} \leq u_n(x_{k_2}), \cdots, M_n^{(k_l)} \leq u_n(x_{k_l})\}$$
$$\to H^*(x_{k_1}, x_{k_2}, \cdots, x_{k_l}), \quad a.s.,$$

*where $H^*(x_{k_1}, x_{k_2}, \cdots, x_{k_l})$ is the marginal distribution of $H(x_1, x_2, \cdots, x_k)$. Especially, for fixed $k \in \mathbb{N}$,*

$$\frac{1}{D_N} \sum_{n=k}^{N} d_n \mathbb{I}\{M_n^{(k)} \leq u_n(x)\} \to H_k(x) = G(x) \sum_{j=0}^{k-1} \frac{(-\log G(x))^j}{j!}, \quad a.s.$$

For bounded Lipschitz 1 functions, we have the following a.s. central limit theorem of order statistics.

**Corollary 2.** *Under the conditions of Theorem 1, for fixed $k \in \mathbb{N}$ and bounded Lipschitz 1 function $f$, we have*

$$\frac{1}{D_N} \sum_{n=k}^{N} d_n f\left(a_n^{-1}\left(M_n^{(k)} - b_n\right)\right) \to \int_{-\infty}^{\infty} f(x) dH_k(x), \quad a.s.$$

# 3 Proofs

As mentioned above, denote levels $u_n(x) = a_n x + b_n, x \in \mathbb{R}, n \geq 1$ and real numbers $x_1 > x_2 > \cdots > x_k$ for fixed $k$. For convenience, let $M_{m,n}^{(j)}$ denote the $j$th maxima of $X_{m+1}, X_{m+2}, \cdots, X_n, 0 \leq m < n$ and $M_n^{(j)} := M_{0,n}^{(j)}$. Set

$$\begin{aligned}\eta_{m,n} &= \mathbb{I}\{M_{m,n}^{(1)} \leq u_n(x_1), M_{m,n}^{(2)} \leq u_n(x_2), \cdots, M_{m,n}^{(k)} \leq u_n(x_k)\} \\ &\quad - P\left(M_{m,n}^{(1)} \leq u_n(x_1), M_{m,n}^{(2)} \leq u_n(x_2), \cdots, M_{m,n}^{(k)} \leq u_n(x_k)\right)\end{aligned}$$

for $0 \leq m < n$ and $\eta_n = \eta_{0,n}$. Before proving the main results, we need some lemmas. Our first lemma provides a bound on the expectation of the difference of the indicator functions for the $j$th maxima of the whole sequence and the $j$th maxima of the subsequence for $j = 1, 2, \cdots, k$.

**Lemma 1.** *Assume (1.1) holds and $m \geq k, n - m \geq k$. Then*

$$\mathbb{E}\left|\mathbb{I}\{M_n^{(j)} \leq u_n(x)\} - \mathbb{I}\{M_{m,n}^{(j)} \leq u_n(x)\}\right| \leq k\frac{m}{n}$$

*uniformly for $j = 1, 2, \cdots, k$ and $x \in \mathbb{R}$.*

**Proof.** First note $\mathbb{I}\{M_n^{(j)} \leq u_n(x)\} - \mathbb{I}\{M_{m,n}^{(j)} \leq u_n(x)\} \neq 0$ if and only if $M_n^{(j)} > M_{m,n}^{(j)}$. The latter implies that $M_m^{(1)} > M_{m,n}^{(j)}$. It is known that the distribution of the general order statistic $M_n^{(j)}$ is

$$P\left(M_n^{(j)} \leq x\right) = \sum_{i=0}^{j-1} \binom{n}{i} (F(x))^{n-i} (1 - F(x))^i,$$

where $\binom{n}{i} = n!/\{i!(n-i)!\}$. Hence,

$$\begin{aligned}&\mathbb{E}\left|\mathbb{I}\{M_n^{(j)} \leq u_n(x)\} - \mathbb{I}\{M_{m,n}^{(j)} \leq u_n(x)\}\right| \\ &\leq P\left(M_n^{(j)} > M_{m,n}^{(j)}\right) \leq P\left(M_m^{(1)} > M_{m,n}^{(j)}\right) \\ &= \sum_{i=0}^{j-1} m\binom{n}{i} \int_{-\infty}^{\infty} (F(x))^{n+m-i-1} (1 - F(x))^i \, dF(x) \\ &\leq \sum_{i=0}^{j-1} m\binom{n}{i} \int_0^1 x^{n+m-i-1}(1-x)^i dx \\ &= \sum_{i=0}^{j-1} m\binom{n}{i} \frac{(n+m-i-1)!i!}{(n+m)!} \\ &\leq j\frac{m}{n} \leq k\frac{m}{n}\end{aligned}$$

uniformly for $1 \leq j \leq k$ and $x \in \mathbb{R}$. □

Our next lemma provides a bound for the covariance of $\eta_m$ and $\eta_n$, which will be used later for estimating the moment of the weighted sum of $\eta_n$.

**Lemma 2.** *Assume (1.1) holds. Then*

$$|\text{Cov}(\eta_m, \eta_n)| \leq 2k^2 \frac{m}{n} \quad (3.1)$$

*for $m \geq k, n - m \geq k$.*

**Proof.** The desired result follows by Lemma 1 and noting

$$\begin{aligned}
|\text{Cov}(\eta_n, \eta_m)| &\leq 2 \left( \mathbb{E} \left| \mathbb{I}\{M_n^{(k)} \leq u_n(x_k)\} - \mathbb{I}\{M_{m,n}^{(k)} \leq u_n(x_k)\} \right| \right. \\
&\quad + \mathbb{E} \left| \mathbb{I}\{M_n^{(k-1)} \leq u_n(x_{k-1})\} - \mathbb{I}\{M_{m,n}^{(k-1)} \leq u_n(x_{k-1})\} \right| \\
&\quad + \cdots \\
&\quad \left. + \mathbb{E} \left| \mathbb{I}\{M_n^{(1)} \leq u_n(x_1)\} - \mathbb{I}\{M_{m,n}^{(1)} \leq u_n(x_1)\} \right| \right). \square
\end{aligned}$$

The following lemma is useful to estimate the moment of the weighted sum of $\eta_n - \eta_{m,n}$.

**Lemma 3.** *Assume (1.1) holds. For $m \geq k, n - m \geq k$, we have*

$$\mathbb{E} |\eta_n - \eta_{m,n}| \leq 2k^2 \frac{m}{n}. \quad (3.2)$$

**Proof.** Note

$$\begin{aligned}
&\mathbb{E} |\eta_n - \eta_{m,n}| \\
&= 2\mathbb{E} \left( \prod_{j=1}^{k} \mathbb{I}\{M_{m,n}^{(j)} \leq u_n(x_j)\} - \prod_{j=1}^{k} \mathbb{I}\{M_n^{(j)} \leq u_n(x_j)\} \right) \\
&\leq 2 \sum_{j=1}^{k} \mathbb{E} \left( \mathbb{I}\{M_{m,n}^{(j)} \leq u_n(x_j)\} - \mathbb{I}\{M_n^{(j)} \leq u_n(x_j)\} \right)
\end{aligned}$$

by the elementary inequality $\left| \prod_{j=1}^{l} y_j - \prod_{j=1}^{l} z_j \right| \leq \sum_{j=1}^{l} |y_j - z_j|$ for all $|y_j| \leq 1, |z_j| \leq 1, j = 1, 2, \cdots, l$. By using Lemma 1, the proof is complete. □

**Lemma 4.** *Under the conditions of Theorem 1, for any $\omega$ with $k \leq m \leq \omega \leq n$ and $p \in \mathbb{N}$,*

$$\mathbb{E}\left|\sum_{l=\omega}^{n} d_l(\eta_l - \eta_{m,l})\right|^p \leq 2^{2p-1}k\left(2 + kp^{\frac{p}{2}}\right)\left(\sum_{l=\omega}^{n} ld_l^2\right)^{\frac{p}{2}}.$$

**Proof.** Note $|\eta_l - \eta_{k,l}| \leq 4$ and, for $m \geq k$, $l - m \geq k$, by using Lemma 3, we have

$$\mathbb{E}|\eta_l - \eta_{m,l}|^p \leq 2 \cdot 4^{p-1}\mathbb{E}|\eta_l - \eta_{m,l}| \leq 4^p k^2\left(\frac{m}{l}\right).$$

Then by Hölder inequality and (1.4), similar to the arguments in Lemma 3 of Hörmann (2006), we have

$$\mathbb{E}\left|\sum_{l=\omega+k}^{n} d_l(\eta_l - \eta_{m,l})\right|^p \leq 4^p k^2 p^{\frac{p}{2}}\left(\sum_{l=\omega}^{n} ld_l^2\right)^{\frac{p}{2}}.$$

By using the $C_r$ inequality,

$$\mathbb{E}\left|\sum_{l=\omega}^{n} d_l(\eta_l - \eta_{m,l})\right|^p$$

$$\leq 2^{p-1}\left(\mathbb{E}\left|\sum_{l=\omega}^{\omega+k-1} d_l(\eta_l - \eta_{m,l})\right|^p + \mathbb{E}\left|\sum_{l=\omega+k}^{n} d_l(\eta_l - \eta_{m,l})\right|^p\right)$$

$$\leq 2^{p-1}\left(2k \cdot 4^p \max_{\omega \leq l \leq \omega+k-1} d_l^p + \mathbb{E}\left|\sum_{l=\omega+k}^{n} d_l(\eta_l - \eta_{m,l})\right|^p\right)$$

$$\leq 2^{2p-1}k\left(2 + kp^{\frac{p}{2}}\right)\left(\sum_{l=\omega}^{n} ld_l^2\right)^{\frac{p}{2}}.$$

The proof is complete. $\square$

The following is the result of Lemma 2, Lemma 4 and slight changes to the proof of Lemma 4 of Hörmann (2006) (or Lemma 2 of Hörmann (2007a)).

**Lemma 5.** *Under the conditions of Theorem 1, for every $p \in \mathbb{N}$, there exists a constant $C_p > 0$ such that*

$$\mathbb{E}\left|\sum_{n=k}^{N} d_n \eta_n\right|^p \leq C_p\left(\sum_{k \leq m \leq n \leq N} d_m d_n \left(\frac{m}{n}\right)^\alpha\right)^{\frac{p}{2}}.$$

The following is the result of Hörmann (2006, 2007a).

**Lemma 6.** *Assume (1.5) holds. For any $\alpha > 0$ and $\eta < \rho$, we have*
$$\sum_{k \leq m \leq n \leq N} d_m d_n \left(\frac{m}{n}\right)^\alpha = O\left(\frac{D_N^2}{(\log D_N)^\eta}\right).$$

**Proof of Theorem 1.** By Lemmas 4 and 5, using Markov inequality and the subsequence procedure, we obtain the desired results, cf. Hörmann (2006, 2007a, 2007b). □